\documentclass[12pt,oneside]{amsart}
\usepackage{latexsym}
\usepackage{amsmath}
\usepackage{amssymb}
\usepackage{ifthen}
\usepackage{amsopn}

\usepackage{epsf}
\usepackage{psfrag}

\newtheorem{theor}{Theorem} 

\newtheorem{theorem}{Theorem}

\newtheorem{thm}{Theorem}[section] \newtheorem{lem}[thm]{Lemma}
 
\newtheorem{prop}[thm]{Proposition} \theoremstyle{definition}
\newtheorem{defn}[thm]{Definition}
\newtheorem{notation}[thm]{Notation}
\newtheorem{conv}[thm]{Convention} \newtheorem{rem}[thm]{Remark}

\begin{document}

\title[Relative hyperbolicity and Artin groups] {Relative
  hyperbolicity and Artin groups} \author[I.~Kapovich]{Ilya Kapovich}

\address{Department of Mathematics, University of Illinois at
  Urbana-Champaign, 1409 West Green Street, Urbana, IL 61801, USA}
\email{kapovich@math.uiuc.edu}

\author[P.~Schupp]{Paul Schupp}

\address{Department of Mathematics, University of Illinois at
  Urbana-Champaign, 1409 West Green Street, Urbana, IL 61801, USA}
\email{schupp@math.uiuc.edu}


%\date{\today}

\begin{abstract}
  Let $G=\langle a_1,\dots, a_n \, |\, a_ia_ja_i\dots = a_ja_ia_j\dots,
  i<j\rangle$ be an Artin group and let $m_{ij}=m_{ji}$ be the length
  of each of the sides of the defining relation involving $a_i$ and
  $a_j$.  We show if all $m_{ij}\ge 7$ then $G$ is relatively
  hyperbolic in the sense of Farb with respect to the collection of
  its two-generator subgroups $\langle a_i, a_j\rangle$ for which
  $m_{ij}<\infty$.
\end{abstract}
\maketitle

\section{Introduction}

The notion of a word-hyperbolic group introduced by Gromov~\cite{Gro}
has played a pivotal role in the development of geometric group theory
for the last fifteen years.  In \cite{Gro,Gro1} Gromov also suggested
a way of generalizing this notion to that of a group \emph{relatively
  hyperbolic} with respect to a collection of subgroups called
\emph{parabolic subgroups}. Several researchers proposed ways of
making these ideas precise. First, Farb~\cite{Fa} defined relative
hyperbolicity by requiring that the Cayley graph of a group with
cosets of the parabolic subgroups ``coned-off'' is a hyperbolic metric
space. Later Bowditch~\cite{Bo97} gave a rigorous interpretation of
Gromov's original approach which mimics the case of a geometrically
finite Kleinian group. Recently Yaman~\cite{Ya} showed that the
Gromov-Bowditch version of relative hyperbolicity can be characterized
in terms of convergence group actions. Bowditch had earlier obtained
such a characterization for word-hyperbolic groups~\cite{Bo98}.
Szczepa\'nski~\cite{Sz} investigated the relationship between
different versions of relative hyperbolicity. He proved that relative
hyperbolicity in the sense of Gromov-Bowditch-Yaman implies that of
Farb but not vice versa.  Bowditch~\cite{Bo97} observed that relative
hyperbolicity in the sense of Farb together with what Farb termed the
``bounded coset penetration property'' implies relative hyperbolicity
in the sense of Gromov-Bowditch-Yaman.  A number of interesting
results regarding relatively hyperbolic groups are obtained in
\cite{Bo99,Bo01,BrFa,Gol,Da,Da1,GL,Hr,MM,Sz1}. In particular,
Szczepa\'nski~\cite{Sz1} recently provided a way of constructing
groups relatively hyperbolic in the sense of Farb by mimicking the
hyperbolization of polyhedra construction.

In this paper we study the class of Artin groups, that is, groups
given by a presentation of the form:

\begin{equation*}\label{art} 
G=\langle a_1,\dots,  a_n \, | u_{ij}=u_{ji}, \text{ for }1\le 
i<j\le n \rangle \tag{\dag} 
\end{equation*} 
where for $i\ne j$
\[ 
u_{ij}:=\underbrace{a_ia_ja_i\dots}_{m_{ij} \text{ terms }}
\] 
and where $m_{ij}=m_{ji}$ for each $i<j$. We allow $m_{ij}=\infty$ in
which case the relation $u_{ij}=u_{ji}$ is omitted from
presentation~\eqref{art}.

In the theory of Artin groups a subgroup of an Artin group $G$
generated by a subset of $\{a_1,\dots, a_n\}$ is called a
\emph{parabolic subgroup} (see, for example, \cite{Par}).  It follows
from the results of Appel-Schupp~\cite{ApSc} that if in \eqref{art}
all $m_{ij}\ge 4$ then for any $i<j$ the two-generator parabolic
subgroup $G_{ij}:=\langle a_i, a_j\rangle\le G$ is itself an Artin
group with the presentation

\[
G_{ij}=\langle a_i, a_j\ |\ u_{ij}=u_{ji}\rangle.
\]

Although we will not use the fact, it is easy to see that for an even
$m_{ij}=2k$ the group $G_{ij}$ is isomorphic to the Baumslag-Solitar
group $B(k,k)=\langle x,y | y^{-1}x^ky=x^k\rangle$. Similarly, if
$m_{ij}=2k+1$ is odd, then $G_{ij}$ is isomorphic to the torus-knot
group $\langle x, y | x^{2k+1}=y^2\rangle$.

An Artin group is said to be of \emph{extra large type} if all $m_{ij}
\ge 4$ in presentation~\eqref{art}.  Even before the general theory of
hyperbolic groups Appel and Schupp \cite{ApSc} proved theorems about
Artin groups of extra large type by showing that they were
``relatively small cancellation'' with respect to the collection of
subgroups $G_{ij}$.  Our main result here is:

\begin{theor}\label{A}
  Let $G$ be an Artin group given by presentation~\eqref{art} above.
  Assume that for all $i<j$ we have $m_{ij}\ge 7$.  Then $G$ is
  relatively hyperbolic in the sense of Farb with respect to the
  collection of subgroups $\{G_{ij} | m_{ij}<\infty\}$, where
  $G_{ij}=\langle a_i, a_j\rangle\le G$.
\end{theor}

Thus Theorem~\ref{A} asserts that $G$ is relatively hyperbolic with
respect to the family of non-free two-generator parabolic subgroups.
This provides additional justification for using the term ``parabolic
subgroup'' in the context of Artin groups.

Since we do not actually use it, we refer to Farb's paper \cite{Fa}
for the precise definition of the ``bounded coset penetration
property.''  But it is easy to see that if two subgroups $H_1,H_2$ of
a group $G$ have infinite intersection then $G$ does not satisfy the
bounded coset penetration property with respect to the collection
$(H_1,H_2)$ (or with respect to any larger collection of subgroups).
Thus for an Artin group $G$ as in Theorem~\ref{A} if there are
distinct $i,j,k$ such that $m_{ij}$ and $m_{ik}$ are finite then $G$
does not have the bounded coset penetration property with respect to
the collection of subgroups $G_{ij}$. This raises the interesting
question of finding different conditions which ensure good
group-theoretic properties in the presence of Farb's definition of
relative hyperbolicty. It seems plausible that most Artin groups
satisfy some sort of a ``nested'' version of the bounded coset
penetration property, especially since many of these groups are known
to be CAT(0) and biautomatic (see for example,
\cite{AC,BC,CD,NR,Peif}).

\section{Artin groups and small cancellation theory}

In this section we explain how one can apply small cancellation theory
to study Artin groups, even though the given finite
presentation~\eqref{art} does not have good small cancellation
properties. In order to do this we assume that the reader is familiar
with the basics of standard small cancellation theory.  For
definitions and terminology see Lyndon and Schupp~\cite{LS}.

Let $R_{ij}$ be the symmetrized set of words in $F(a_i,a_j)$ generated
by the defining relator, $u_{ij}u_{ji}^{-1}$, and its inverse.  The
cancellation between two noninverse elements of $R_{ij}$ can be almost
half of each of the words but \emph{not} an entire half. Consider, for
example, a product such as

\[ (a_1 a_2 a_1 a_2 a_1^{-1} a_2^{-1} a_1^{-1} a_2^{-1}) (a_2 a_1 a_2 a_1^{-1} 
a_2^{-1} a_1^{-1} a_2^{-1} a_1) \]
   
It is not difficult to verify that the set $R_{ij}$ satisfies the
``flat'' small cancellation condition $C(4)-T(4)$. To exploit this
fact we must use the right normal form for elements of $F(a_1,\dots,
a_n)$.

\begin{conv}
  
  For the remainder of this paper we denote $F:=F(a_1,...,a_n)$.  In
  the free group $F$ every nontrivial reduced word $w$ has a unique
  \emph{normal form with exponents}
 
\[ w = a^{n_1}_{h_1} ... a^{n_s}_{h_s} \] 
 
\noindent where each $h_t \ne h_{t+1} $ and each $n_t \ne 0$. The integer $s$ is
the \emph{syllable length} of $w$ and is written $ ||w||$. The
subwords $a_{h_i}^{n_i}$ are called the \emph{syllables} of $w$.
\end{conv}

\begin{notation} 
  For each $i<j$ such that $m_{ij}<\infty$ let $\mathcal{R}_{ij}$ be the
  set of \emph{all} nontrivial cyclically reduced words in
  $F(a_i,a_j)$ which are equal to the identity in the group $G_{ij}$.
\end{notation}

By introducing ``strips'' to study the fine geometry of C(4)-T(4)
diagrams, Appel and Schupp~\cite{ApSc} proved:

\begin{prop}\label{prop:AS} 
  Suppose $2\le m_{ij}<\infty$ and let $w$ be a word from
  $\mathcal{R}_{ij}$. Then $||w||\ge 2m_{ij}$.
\end{prop}

\begin{notation} 
  To deal with the Artin group $G$ we now switch to the following
  \emph{infinite} presentation of $G$:
 
\begin{equation*}\label{art1} 
G=\langle a_1,\dots,  a_n \, |\, {\mathcal R} \rangle. \tag{\ddag} 
\end{equation*} 
where
\[ 
{\mathcal R} =\underset{m_{ij}<\infty} \bigcup \mathcal{R}_{ij}.
\] 
\end{notation} 
  
\begin{rem}
  The point of shifting to the infinite presentation \eqref{art1} is
  that it allows a strong use of minimality. In considering van Kampen
  diagrams for a word $w = 1$ in $G$, it suffices to consider
  \emph{minimal} diagrams, that is, diagrams with as few as regions as
  possible.  In a minimal diagram $\Delta$ over~\eqref{art1} distinct
  regions labelled by relators from the same $\mathcal{R}_{ij}$ cannot
  have even a vertex in common since they could then be combined into
  a single region, contradicting minimality. Also, any common edge
  shared by regions labelled from some $\mathcal{R}_{ij}$ and
  $\mathcal{R}_{il}$ where $j \ne l$ must be labelled by a power of
  the generator $a_i$, that is, the edge label has syllable length
  one.  If $G$ is of extra large type, Proposition~\ref{prop:AS} will
  show that minimal diagrams for the presentation~\eqref{art1} have
  ``hyperbolic'' $C(8)$- geometry.
\end{rem}

\begin{notation}
  As in \cite{LS}, if $\Delta$ is a \emph{map} and $D$ is a region of
  $\Delta$, then $i(D)$ will denote the \emph{interior degree} of $D$,
  that is, the number of $\Delta$-interior edges of $D$, where each
  edge occurring twice in a boundary cycle of $D$ is counted twice in
  $i(D)$. Recall also that a region $D$ of $\Delta$ is said to be a
  \emph{boundary region} if $\partial D\cap \partial \Delta \ne
  \emptyset$. Note that a boundary region need not contain a boundary
  edge of $\Delta$ as the intersection may contain only vertices of
  $\partial \Delta$. We will say that a boundary region $D$ is a
  \emph{simple boundary region} of $\Delta$ if the intersection
  $\partial D\cap \partial \Delta\ne \emptyset$ is connected and if
  the edges of $\partial D\cap \partial \Delta$ are consecutively
  traversed in some boundary cycle of $\Delta$.
\end{notation}
Lemma~4.1 in Ch.~V of \cite{LS} asserts that in a $(6,3)$-map the
boundary of every region is a simple closed curve. In particular this
means that if $D$ is a simple boundary region of a $(6,3)$-map
$\Delta$ then the intersection $\partial D\cap \partial \Delta$ is
consecutively traversed in some boundary cycle of $D$.

We shall need the following ``basic fact'' of small cancellation
theory~\cite{LS}:

\begin{prop}\label{use}
  
  Let $\Delta$ be a $(6,3)$-map without vertices of degree one and
  with no simple boundary regions of interior degree zero.  Then at
  least one of the following occurs:

\begin{enumerate}
\item[(a)] The map $\Delta$ contains at least two simple boundary
  regions of interior degree one.
  
\item[(b)] The map $\Delta$ has at least three simple boundary regions
  of interior degree at most three.
\end{enumerate}
\end{prop}

Proposition~\ref{prop:AS} and Proposition~\ref{use} together
imply~\cite{ApSc}:

\begin{prop}\label{AS}
  Suppose $G$ is given by \eqref{art} where all $m_{ij}\ge 4$.  Then:
\begin{enumerate}
\item Every minimal diagram over the infinite presentation
  \eqref{art1} satisfies $C(8)$ and every interior edge in such a
  diagram is labeled by a power of some generator $a_i$.
\item If $w$ is a nontrivial freely reduced word representing $1$ in
  $G$ then $w$ contains a subword $v$ such that $v$ is also a subword
  of some $r \in \mathcal{R}_{ij}$ with $r =vu$, $||u||\le 3$ and
  $||v||\ge 2m_{ij}-3$.
\item If $m_{ij}, m_{ik}<\infty$ for some $k\ne j$ then $G_{ij}\cap
  G_{ik}=\langle a_i \rangle$.  Moreover, if $\{i,j\}\cap \{t,
  k\}=\emptyset$ then $G_{ij}\cap G_{tk}=\{1\}$ in $G$.
\end{enumerate} 
\end{prop}

\begin{conv} For the remainder of this section we assume that $G$ is given by presentation~\eqref{art} where all $m_{ij}\ge 4$.
\end{conv}

\begin{defn} Let $w$ be a freely reduced word in $F(a_1,\dots, a_n)$.
  We say that $w$ is \emph{Artin-reduced} if $w$ does not contain a
  subword $v$ such that $v$ is also a subword of some
  $r\in\mathcal{R}_{ij}$ with $r =vu$, $||u||\le 3$.
  
  Similarly, a word $w$ is \emph{strongly Artin-reduced} if $w$ does not
  contain a subword $v$ such that $v$ is also a subword of some
  $r\in\mathcal{R}_{ij}$ with $r =vu$, $||u||\le 4$.
\end{defn}

An equality diagram for a pair of
Artin-reduced words has the standard structure, consisting of
``islands'' in which each region ``cuts through'' the diagram from
``top to bottom'' as described precisely in the following lemma.

\begin{figure}[here]

  \psfrag{u}{$u$} \psfrag{v}{$v$}
  
  \epsfxsize=4.5in \centerline{\epsfbox{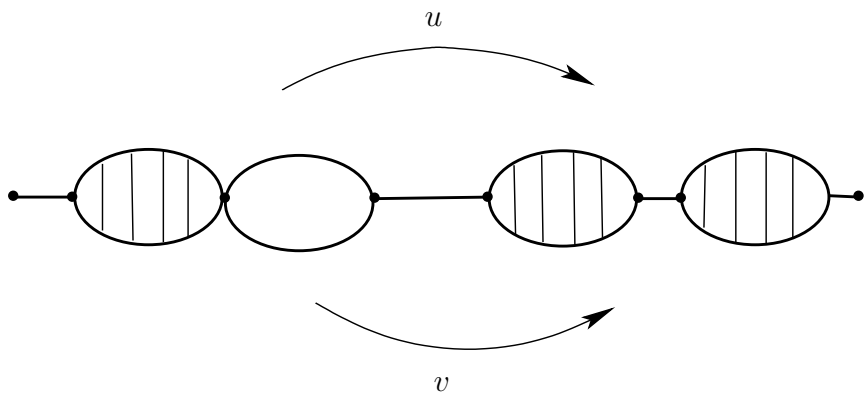}} \caption{Equality
    diagram}\label{Fi:eq2}
\end{figure}

\begin{figure}[here]

  \psfrag{u'}{$u'$} \psfrag{v'}{$v'$}
  
  \psfrag{Q1}{$Q_1$}
  
  \psfrag{Q2}{$Q_2$} \psfrag{Qt-1}{$Q_{t-1}$} \psfrag{Qt}{$Q_t$}
  \psfrag{p1}{$p_1$} \psfrag{p2}{$p_2$} \psfrag{pt-2}{$p_{t-2}$}
  \psfrag{pt-1}{$p_{t-1}$} \epsfxsize=4.5in
  \centerline{\epsfbox{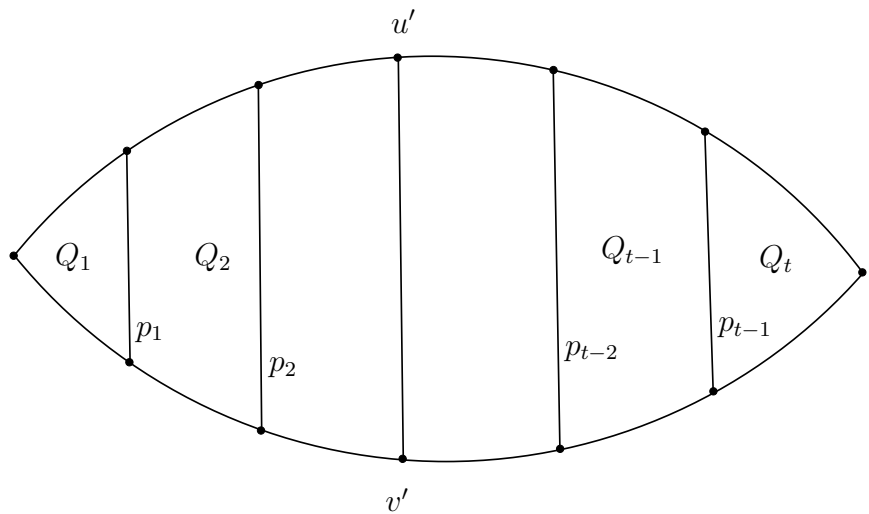}} \caption{A disk component in an
    equality diagram}\label{Fi:eq}
\end{figure}

\begin{lem}\label{shape}
  Suppose $\Delta$ is a minimal diagram over \eqref{art1} such that
  the boundary cycle of $\Delta$ corresponds to the relation $u=_G v$
  where both $u$ and $v$ are Artin-reduced. Then $\Delta$ is a union
  of several (possibly none) linearly arranged disk components
  connected by (possibly degenerate) arcs, as shown in
  Figure~\ref{Fi:eq2}. Moreover, each of these disk components
  corresponds to an equality $u'=_G v'$, where $u'$ is a subword of
  $u$ and $v'$ is a subword of $v$, and the component has the form
  shown in Figure~\ref{Fi:eq} satisfying the following properties:
\begin{enumerate}
  
\item Each region $Q_s$ corresponds to a relator from \eqref{art1}.
\item Each interior edge $p_s$ is labeled by a single syllable
  $a_{i_s}^{n_s}$ and has one endpoint on the ``top'' part of
  $\partial \Delta$ corresponding to $u$ and its other endpoint on the
  ``bottom'' part of $\partial \Delta$ corresponding to $v$.
\item Each region $Q_s$ has at least one boundary edge contained in
  the top part of $\partial \Delta$ and at least one boundary edge
  contained in the bottom part of $\partial \Delta$.
\end{enumerate}
\end{lem}

\begin{proof}

  By removing interior vertices of degree two and combining edge
  labels, we follow the usual convention that $\Delta$ has no interior
  vertices of degree two.  Since $u$ and $v$ are Artin-reduced, it is
  clear that $\Delta$ is indeed a union of linearly arranged disk
  components, each corresponding to an equality diagram for $u'=_G v'$
  for some subwords $u'$ of $u$ and $v'$ of $v$, as shown in
  Figure~\ref{Fi:eq2}. Also, each of the disk components is a minimal
  diagram over \eqref{art1}.

  Recall that by Proposition~\ref{AS} in a minimal diagram over
  \eqref{art1} no two regions corresponding to relators from the same
  $\mathcal{R}_{ij}$ have a common edge.  Moreover, if $r\in
  \mathcal{R}_{ij}, r'\in \mathcal{R}_{il}$, where $l\ne j$, have a
  nontrivial common subword, this subword must be a power of $a_i$.
  Thus every interior edge of a minimal diagram over \eqref{art1} is
  labeled by some $a_i^s$.

\begin{figure}[here]
  
  \psfrag{u'}{$u'$} \psfrag{v'}{$v'$} \psfrag{x1}{$x_1$}
  \psfrag{x2}{$x_2$} \psfrag{Q}{$Q$} \psfrag{Q'}{$Q'$} \psfrag{S}{$S$}
  \psfrag{S'}{$S'$} \psfrag{q}{$q$} \psfrag{q'}{$q'$}\psfrag{R}{$R$}
  \psfrag{R'}{$R'$}
  
  \epsfxsize=4.5in \centerline{\epsfbox{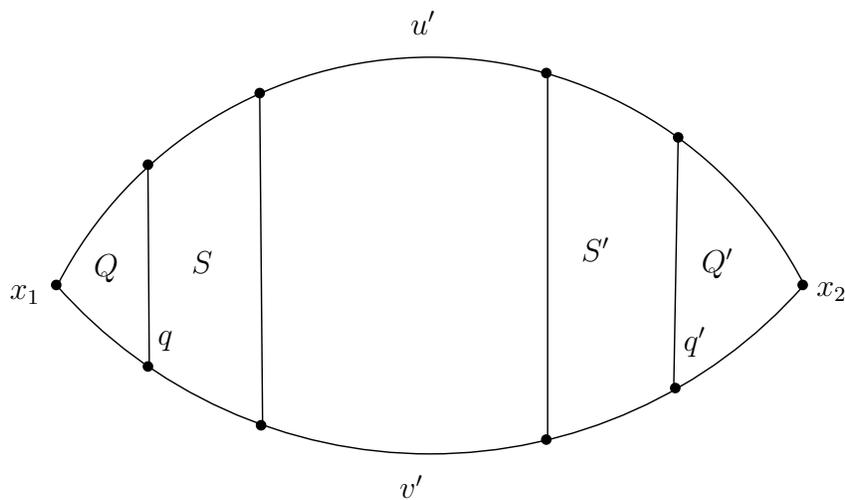}}

\caption{Equality diagram $\Delta_0$ for $u'=_G v'$ in Lemma~\ref{shape}}\label{Fi:eq1}
\end{figure}

Consider an individual disk component $\Delta_0$ corresponding to an
equality $u'=_G v'$, where $u'$ is a subword of $u$ and $v'$ is a
subword of $v$.

Let $x_1$ be the common initial vertex of the paths labeled $u'$ and
$v'$ in the boundary of $\Delta_0$. Similarly, let $x_2$ be their
common terminal vertex. If $\Delta_0$ has only one region,
Lemma~\ref{shape} certainly holds. Thus we now suppose that $\Delta_0$
has at least two regions. Since $\Delta_0$ is minimal over
presentation~\eqref{art1}, by Proposition~\ref{AS} $\Delta_0$ is a
$C(8)$-diagram. Hence we can apply part Proposition~\ref{use} to
$\Delta_0$. Since $u'$ and $v'$ are Artin-reduced, case (b) of
Proposition~\ref{use} cannot occur. Thus case (a) takes place and
there are two distinct boundary regions $Q$ and $Q'$ in $\Delta_0$
each with interior degree one and with a single syllable written on
its interior edge, as shown in Figure~\ref{Fi:eq1}. Denote those edges
by $q$ and $q'$ accordingly.  It is clear that $x_1$ belongs to one of
these regions and $x_2$ belongs to the other. Suppose $x_1$ is in $Q$
and $x_2$ is in $Q'$. If $\Delta_0$ consists of two or three regions,
the statement of the lemma obviously holds. Otherwise let $\Delta'$ be
obtained from $\Delta_0$ by first removing the regions $Q,Q'$, and, if
their removal creates any vertices of degree one, then successively
removing vertices of degree one and their incident edges until no such
vertices remain. If $\Delta'$ has fewer than two regions then the
statement of the lemma holds.

Otherwise the diagram $\Delta'$ contains at least two regions and we
can apply Proposition~\ref{use} to it. Again, clearly case (b) of
Proposition~\ref{use} cannot occur and case (a) takes place. Let
$S,S'$ be two distinct boundary regions in $\Delta'$, each of interior
degree one and thus with a single syllable written on the interior
edge. Since $u$ and $v$ are weakly Artin-reduced, one of these
regions, say $S$, contains the edge $q$ and the other region $S'$
contains $q'$.
  
We claim that each of $S,S'$ has an edge in common with both the top
and the bottom portions of the boundary of $\Delta_0$. If this is not
the case then either $u'$ or $v'$ contains a subword of $r_{ij}^*$
missing at most two letters, contradicting the assumption that $u'$
and $v'$ are Artin-reduced. Hence we can remove $S$ and $S'$ from
$\Delta'$ and repeatedly apply the same argument to the remaining
diagram.
  
This yields the statement of Lemma~\ref{shape}.
\end{proof}  

\section{Coned-off Cayley graphs and relative hyperbolicity}

Recall that a \emph{geodesic path} in a metric space $(X,d)$ is an
isometric embedding $\gamma: [a,b]\to X$, where $[a,b]\subseteq
\mathbb R$. In this situation we call the set $\gamma([a,b])$ a
\emph{geodesic segment} from $x=\gamma(a)$ to $y=\gamma(b)$ in $X$ and
denote it by $[x,y]$. A metric space $(X,d)$ is said to be
\emph{geodesic} if for any $x,y\in X$ there exists a geodesic segment
$[x,y]$ in $X$.  The general notion of a ``hyperbolic metric space''
is, of course, due to Gromov~\cite{Gro}:

\begin{defn}[Hyperbolic Metric Space]
  A geodesic metric space $(X,d)$ is said to be \emph{hyperbolic} if
  there is a number $\delta\ge 0$ such that for any three points
  $x,y,z\in X$, for any geodesic segments $[x,y]$, $[x,z]$ and $[y,z]$
  and for any point $p\in [x,y]$ there is a point $q\in [x,z]\cup
  [y,z]$ such that
\[
d(p,q)\le \delta.
\]
\end{defn}

We refer the reader to \cite{ABC,CDP,GH,Gro} for the background
material on hyperbolic spaces and hyperbolic groups.

We need to recall Farb's definitions of coned-off Cayley graphs and
relative hyperbolicity~\cite{Fa}:

\begin{defn}[Coned-off Cayley Graph]\cite{Fa}
  Let $G$ be a group with a finite generating set $A$ and let
  $H_1,\dots, H_t$ be a family of subgroups of $G$. Let
  $\Gamma=\Gamma_A(G)$ be the Cayley graph of $G$ with respect to $A$.
  
  The \emph{coned-off Cayley graph} $X =\widehat \Gamma_A (G;
  H_1,\dots, H_t)$ is obtained from $\Gamma$ as follows: for each
  coset $gH_i$ (where $g\in G$, $1\le i\le t$) we add a new vertex
  $v(gH_i)$ to $\Gamma$ and for every element $g'\in gH_i$ add an edge
  of length $1/2$ from $v(gH_i)$ to $g'$.
\end{defn}

The vertices $v(gH_i)$ are referred to as \emph{cone vertices} and the
edges from $v(gH_i)$ to elements of $gH_i$ are called \emph{cone
  edges}.

\begin{defn}[Relative Hyperbolicity]\cite{Fa}
  Let $G$ be a finitely generated group and let $H_1,\dots, H_t$ be a
  family of subgroups of $G$. We say that $G$ is \emph{relatively
    hyperbolic} in the sense of Farb with respect to the $t$-tuple
  of subgroups $(H_1,\dots, H_t)$ if for some finite generating set
  $A$ of $G$ the coned-off Cayley graph $X = \widehat \Gamma_A (G;
  H_1,\dots, H_t)$ is a hyperbolic metric space.
\end{defn}

Recall that a finitely generated group $G$ is \emph{hyperbolic} if for
some (and hence for any) finite generating set $A$ of $G$ the Cayley
graph $\Gamma(G,A)$ is a hyperbolic metric space. Thus $G$ is
hyperbolic if and only if it is relatively hyperbolic with respect to
the trivial subgroup $H=\{1\}$.

It is shown in \cite{Fa} that $G$ being relatively hyperbolic with
respect to the given collection of subgroups does not depend on the
choice of the generating set $A$.

\section{Geodesics in the coned-off Cayley graph of an Artin group}

Let $G$ be given by presentation~\eqref{art} where $m_{ij}\ge 7$ for
all $i\ne j$. Let $\Gamma=\Gamma_A(G)$ where $A=\{a_1,\dots, a_n\}$.
Let $X=\widehat \Gamma_A(G; \{G_{ij}| m_{ij}<\infty\})$ be the
coned-off Cayley graph of $G$. If $e$ is an edge of a graph, we will
denote by $o(e)$ the initial vertex of $e$ and by $t(e)$ the terminal
vertex of $e$.  These notations and conventions will be fixed for the
remainder of the paper.

\begin{conv}
  Let $g\in G$ be an arbitrary element and let $\alpha$ be a geodesic
  in $X$ from $1\in G$ to $g$. The goal of this section is to show
  that $\alpha$ is 4-close in $X$ to a path $\gamma$ from $1$ to $g$
  in $\Gamma$ whose label is strongly Artin-reduced. As an
  intermediate step in constructing this path $\gamma$ we will first
  define an auxiliary path $\beta$.  Unless specified otherwise, we
  shall fix $g$ and $\alpha$ for the remainder of this section.
\end{conv}

\subsection{Construction of the path $\beta$}

The path $\alpha$ can be subdivided as a concatenation of paths
\[
\alpha=\alpha_1\dots \alpha_t
\]
where $t\ge 1$ and where each $\alpha_k$ is either a path in $\Gamma$
from a vertex of $\Gamma$ to a vertex of $\Gamma$, in which case
$\alpha_k$ is referred to as a \emph{$\Gamma$-block}, or $\alpha_k=e_k
e_k'$, where $e_k$ and $e_k'$ are cone-edges separated by a
cone-vertex, in which case $\alpha_k$ is referred to as a
\emph{cone-block}.  Thus $t(e_k)=o(e_k')$ is a cone-vertex of $X$ and
$t(e_k')=o(e_k)g_{k}$ for some $g_{k}\in G_{ij}-\{1\}$. Moreover, we
can assume that the $\Gamma$-blocks are chosen maximally, so that no
two such blocks are consecutive in $\alpha$.  For each cone-block
$\alpha_k$ choose a freely reduced word $v_k$ in $G_{ij}$ of the
smallest possible syllable length representing $g_k$.  Note that if
$g_k$ belongs to two distinct $G_{ij}$-subgroups and thus is a power
of the generator then the pair $\{i,j\}$ may not be uniquely defined
by $g_k$, but in that case the word $v_k$ of syllable length one will
be uniquely determined by $g_k$.

\begin{rem}
  The choice of the word $v_k$ for the cone-block $\alpha_k$ as a word
  in the corresponding $G_{ij}$ which is ``geodesic'' with respect to
  syllable length is a crucial feature of our argument.
\end{rem}

If $\alpha_k$ is a $\Gamma$-block then the label $v_k$ of $\alpha_k$
is an $A$-geodesic word.

\begin{notation}
  The word
\[
v=v_1\dots v_t
\]
certainly represents the element $g$.  We denote the path in $\Gamma$
from $1$ to $g$ with label $v$ by $\beta$.  Thus
\[
\beta=\beta_1\dots \beta_t
\]
where $\beta_k$ is the path in $\Gamma$ with label $v_k$ from the
initial point of $\alpha_k$ to the terminal point of $\alpha_k$ .
\end{notation}

The following lemma is a straightforward consequence of the
definitions:

\begin{lem}\label{list}
  Let $\alpha$ and $\beta$ be as above. Then:

\begin{enumerate}
\item The label $v_k$ of each $\Gamma$-block $\alpha_k$ is a
  $\Gamma$-geodesic word.
\item If $v'$ is a subword of the label $v_k$ of some $\Gamma$-block
  such that $v'$ is a word in some $G_{ij}$ then $|v'|\le 1$.
\item Suppose that $\alpha_{k+1}$ is a $\Gamma$-block and that
  $\alpha_k$ is a cone-block with $v_k\in G_{ij}$.  Then the first
  syllable of $v_{k+1}$ is not a power of either $a_i$ or $a_j$.
\item Suppose that $\alpha_k$ is a $\Gamma$-block and that
  $\alpha_{k+1}$ is a cone-block with $v_{k+1}\in G_{ij}$.  Then the
  last syllable of $v_{k}$ is not a power of either $a_i$ or $a_j$.
\item~\label{list:diff} If $\alpha_k, \alpha_{k+1}$ are cone-blocks
  and $v_{k}\in G_{ij}$, $v_{k+1}\in G_{st}$ then $\{i,j\}\ne \{s,
  t\}$.
\item \label{list:close} The paths $\alpha$ and $\beta$ are
  2-Hausdorff close in $X$.
\end{enumerate}
\end{lem}

\begin{lem}

  Suppose that $\alpha_{k-1},\alpha_{k}, \alpha_{k+1}$ are cone-blocks
  with $v_{k-1}\in G_{ij}$, $v_k\in G_{js}$, $v_{k+1}\in G_{sq}$. If
  the first syllable of $v_k$ is a power of $a_j$ and the last
  syllable of $v_k$ is a power of $a_s$ then $||v_k||\ge 3$.
\end{lem}
\begin{proof}
  
  If $||v_k||=2$ then we can replace the path
  $\alpha_{k-1}\alpha_k\alpha_{k+1}$ of length 3 in $X$ by a path of
  length 2 in $X$ consisting of two pairs of cone-edges: the first
  corresponding to $G_{i,j}$ and the second corresponding to $G_{sq}$.
  This contradicts the assumption that $\alpha$ is a geodesic in $X$.
\end{proof}

\subsection{Construction of the path $\gamma$}

Note that the word $v$ is not necessarily reduced. Indeed, it is
possible that $\alpha_k, \alpha_{k+1}$ are two consecutive cone blocks
such that the last syllable of $v_k$ and the first syllable of
$v_{k+1}$ are powers of the same $a_i$. (The above lemmas imply that
this is the only way in which $v$ may fail to be freely reduced).  We
will modify $v$ and $\beta$ to rectify this problem by ``condensing
syllables'' from left to right as follows.

We define a sequence of words $u_1,\dots, u_t$ inductively.

Put $u_1:= v_1$.  Suppose $1\le k<t$ and $u_1,\dots, u_{k}$ are
already defined.  We have two cases to consider:

\noindent (a) Suppose $\alpha_k$ is a cone-block with $v_k\in G_{ij}$.

Let $u_k = z_k s_k$ where $s_k$ is the last syllable of $u_k$.

If $\alpha_{k+1}$ is also a cone-block, let $v_{k+1} = s_{k+1}y_{k+1}$
where $s_{k+1}$ is the first syllable of $v_{k+1}$.  If the last
syllable of $v_{k}$ and the first syllable of $v_{k+1}$ are not powers
of the same generator, then let $u_{k+1} = v_{k+1}$.  If the last
syllable of $v_{k}$ and the first syllable of $v_{k+1}$ are powers of
the same generator, then let $u_k '$ be the reduced form of $u_k
s_{k+1}$ and redefine $u_k$ to be $u_k'$.  Set $u_{k+1}:= y_{k+1}$.

If $\alpha_{k+1}$ is a $\Gamma$-block, put $u_{k+1}:=v_{k+1}$.

\noindent (b) If $\alpha_{k}$ is a $\Gamma$-block, put $u_{k+1}:=v_{k+1}$.

\begin{notation}
  We set
\[
u:=u_1\dots u_t, \tag{$\spadesuit$}
\]
and denote by $\gamma$ the path in $\Gamma$ from $1$ to $g$ with label
$u$.  Also, denote by $\gamma_k$ the subpath of $\gamma$ corresponding
to $u_k$, so that
\[
\gamma=\gamma_1 \dots \gamma_t.
\]
\end{notation}

We summarize the relevant properties of $u$ and $\gamma$:

\begin{lem}\label{gamma} The following hold:
\begin{enumerate}
\item The word $u$ as in ($\spadesuit$) is freely reduced.
\item\label{short} Suppose $\alpha_k$ is a cone-block with $u_k\in
  G_{ij}$. Then for any word $w$ in $G_{ij}$ representing the same
  element as $u_k$ we have $||u_k||\le ||w||+1$.
\item For any $k<t$ we have $||u_ku_{k+1}||=||u_k||+||u_{k+1}||$.
\item For any $\Gamma$-block $\alpha_k$ we have
  $\gamma_k=\beta_k=\alpha_k$ and $u_k=v_k$.
\item\label{close} The paths $\gamma$ and $\beta$ are 2-Hausdorff
  close in $X$.
\end{enumerate}
\end{lem}

\begin{lem}\label{good}
  Let $w$ be a reduced word in some $G_{ij}$, $m_{ij}<\infty$ such
  that $||w||\ge 4$. Suppose $w$ is a subword of $u$ and let $w'$ be
  obtained from $w$ by removing the first and the last syllables of
  $w$. Then $w'$ is a subword of some $u_k$ such that $\alpha_k$ is a
  cone-block of $\alpha$.
\end{lem}

\begin{proof}
  Suppose first that the occurrence of $w$ in $u$ overlaps some
  $u_k=v_k$ corresponding to a $\Gamma$-block $\alpha_k$ of $\alpha$.
  Then some syllable $s$ of $w$ is a subword of this $u_k$.  Recall
  that no two $\Gamma$-blocks are adjacent in $\alpha$. Part~2 of
  Lemma~\ref{list} implies that neither of the syllables preceding or
  following $s$ of $w$ can overlap $u_k$. Let $u_k=asb$. Among the
  words $a,b$ choose the one of largest syllable length. By symmetry
  we may assume that $||a||\ge ||b||$. Hence $||a||\ge 2$ since
  $||w||\ge 4$.  Since $||a||>0$ then $\alpha_{k-1}$ is a cone-block.
  Also, since $||a||\ge 2$ and $w$ is a word in $G_{ij}$, we conclude
  that both $s$ (which is the first syllable of $u_k$) and $u_{k-1}$
  are words in $G_{ij}$ (this is true even if $||u_{k-1}||=1$). This
  is impossible by Lemma~\ref{list}.
  
  Thus $w$ does not overlap any $u_k$ corresponding to $\Gamma$-blocks
  $\alpha_k$.

  Assume now that there is some $u_k$ contained in $w$ such that $w$
  also overlaps $u_{k-1}$ and $u_{k+1}$. Each of $\alpha_{k-1}$,
  $\alpha_{k}$ and $\alpha_{k+1}$ is a cone-block. If $||u_k||=1$ then
  either the part of $w$ following $u_k$ or the part of $w$ preceding
  $u_k$ has syllable length at least two. Since $w$ is a word in
  $G_{ij}$ this yield a contradiction with part~\ref{list:diff} of
  Lemma~\ref{list} and the definition of $u$. Hence $||u_{k}||\ge 2$.
  If either the either the part of $w$ following $u_k$ or the part of
  $w$ preceding $u_k$ has syllable length at least two, then either
  both $v_{k-1}, v_{k}$ or both $v_{k}, v_{k+1}$ belong to $G_{ij}$.
  Again, this is impossible by part~\ref{list:diff} of
  Lemma~\ref{list}. The statement of Lemma~\ref{good} now follows.

  Suppose now that there is no $u_k$ contained in $w$ such that $w$
  also overlaps $u_{k-1}$ and $u_{k+1}$. Thus $w$ overlaps at most two
  of the words $u_k$. If $w$ is contained in a single $u_k$, the
  statement of Lemma~\ref{good} obviously holds. Assume now that $w$
  is a subword of $u_{k} u_{k+1}$ and that $w$ overlaps both of these
  words. If both overlaps have syllable length at least two, then
  $g_k, g_{k+1}\in G_{ij}$, contrary to part~\ref{list:diff} of
  Lemma~\ref{list}. If one of the overlaps has syllable length one,
  then all but the first or the last syllable of $w$ is contained in
  either $u_k$ or $u_{k+1}$ and the statement of Lemma~\ref{good}
  holds.
\end{proof}

\begin{prop}\label{tech}
  The word $u$ is strongly Artin reduced and the paths $\alpha$ and
  $\gamma$ are $4$-Hausdorff close in $X$.
\end{prop}

\begin{proof}
  The paths $\alpha$ and $\gamma$ are $4$-close since each is
  $2$-close to the path $\beta$ by part~\ref{list:close} of
  Lemma~\ref{list} and by part~\ref{close} of Lemma~\ref{gamma}.

  We now show that $u$ is strongly Artin reduced.  If not, there is a
  nontrivial relator $r\in \mathcal{R}_{ij}$ for some $m_{ij}<\infty$
  and a subword $w$ of $u$ such that $r=wy$ and $||y||\le 4$. Thus $w=
  y^{-1}$ in $G$.

  It is now that we use the assumption that all $m_{ij}\ge 7$.  By
  Proposition~\ref{prop:AS} this condition implies that $||r||\ge 14$
  and hence $||w||\ge 10$.  Let $w = s_1 w' s_2$ where $s_1,s_2$ are
  respectively the first and the last syllables of $w$.  Then
  $||w'||\ge 8$ and by Lemma~\ref{good} $w'$ is a subword of some
  $u_k$ such that $\alpha_k$ is a cone-block of $\alpha$.
  
  Moreover, if we write $u_k$ as $u_k=zw'z'$ then
  $||u_k||=||z||+||w'||+||z'||$.  Using the relator $r$, we have
 \[w' =_G s_1^{-1} y^{-1} s_2^{-1} \text{ and hence } u_k =_G z  s_1^{-1} y^{-1} s_2^{-1} z' \]
 But
\[ ||z  s_1^{-1} y^{-1} s_2^{-1} z' || \le ||z|| + 1 + 4 + 1 + ||z'||<||z||+||w'|| -1 +||z'||=||u_k|| -1,
 \]
 where the last inequality holds since $||w'||\ge 8$.  Hence there
 exists a word representing the same element as $u_k$ but with
 syllable length less than $||u_k||-1$.
 
 This contradicts part~\ref{short} of Lemma~\ref{gamma}.
\end{proof}

\section{Proof of the main result}

We recall the following useful fact due to Papasoglu~\cite{Pa}:

\begin{prop}[Thin bigons criterion]\label{pap}
  Let $C$ be a connected graph with the simplicial metric $d$.  Then
  $C$ is hyperbolic if and only if there is some $\delta>0$ such that
  for any points $x,y\in C$ (possibly inside edges) any two geodesic
  paths from $x$ to $y$ in $C$ are $\delta$-Hausdorff close.
\end{prop}

\begin{rem}
  The above result was stated in \cite{Pa} only for Cayley graphs of
  finitely generated groups. However, it is easy to see that
  Papasoglu's proof~\cite{Pa} does not use the Cayley graph assumption
  and works for any connected graph with the simplicial metric. This
  was noted, for example, in~\cite{NS}.
\end{rem}

We are now ready to prove Theorem~\ref{A}.

\begin{theorem}\label{thm:rel}
  Let $G$ be an Artin group given by presentation \eqref{art}.
  Suppose that for all $i<j$ we have $m_{ij}\ge 7$.
  
  Then $G$ is relatively hyperbolic in the sense of Farb with respect
  to the collection of subgroups $\{G_{ij} | m_{ij}<\infty\}$.
\end{theorem}

\begin{proof}
  Put $A=\{a_1,\dots, a_n\}$.  Let $\Gamma=\Gamma_A(G)$ and let
  \[X=\widehat \Gamma_A (G;\{G_{ij}| m_{ij}<\infty\})\] be the
  coned-off Cayley graph of $G$. We need to show that $X$ is
  hyperbolic.
  
  Note that if we barycentrically subdivide each of the $\Gamma$-edges
  of $X$ to obtain the graph $X'$, then twice the metric of $X$
  coincides with the simplicial metric for $X'$. Hence to establish
  Theorem~\ref{thm:rel} by using Papasoglu's criterion it suffices to
  show that geodesic bigons in $X$ (with endpoints possibly inside
  edges) are uniformly thin.  Moreover, the definition of $X$ implies
  that it is enough to prove the following:
  
  {\bf Claim.} There exists a constant $K>0$ with the following
  property. Let $\alpha_1, \alpha_2$ be geodesics in $X$ from $h_1\in
  G$ to $g_1\in G$ and from $h_2\in G$ to $g_2\in G$ respectively such
  that:
\begin{enumerate}
\item Either $h_1^{-1}h_2\in G_{ij}$ for some $m_{ij} <\infty$ or
  $d_\Gamma(h_1,h_2)\le 1$.
\item Either $g_2^{-1}g_1\in G_{sl}$ for some $m_{sl} <\infty$ or
  $d_\Gamma(g_1,g_2)\le 1$.
\end{enumerate}
Then $\alpha_1$ and $\alpha_2$ are $K$-Hausdorff close in $w$.

\begin{figure}[here]

  \psfrag{a}{$a$} \psfrag{b}{$b$} \psfrag{h1}{$h_1$}
  \psfrag{h2}{$h_2$}
  
  \psfrag{g1}{$g_1$}
  
  \psfrag{g2}{$g_2$} \psfrag{alpha1}{$\alpha_1$}
  \psfrag{alpha2}{$\alpha_2$} \psfrag{gamma1}{$\gamma_1$}
  \psfrag{gamma2}{$\gamma_2$} \psfrag{lambda1}{$\lambda_1$}
  \psfrag{lambda2}{$\lambda_2$} \epsfxsize=4.5in
  \centerline{\epsfbox{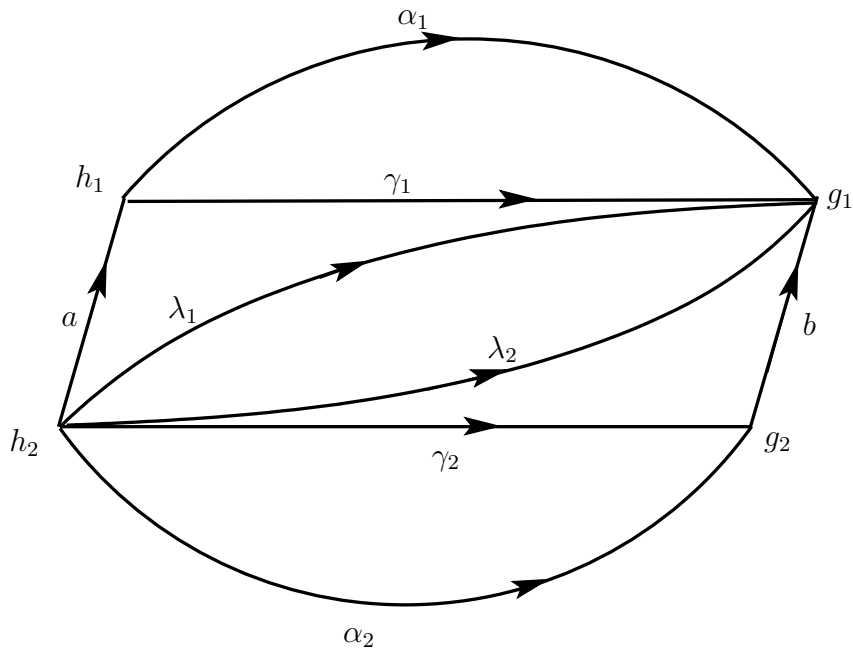}}
\caption{Thinness of bigons}\label{Fi:main}
\end{figure}

To verify the Claim, let $a, b\in G$ be such that $h_1=h_2a$ and
$g_1=g_2b$.  By Proposition~\ref{tech} for $l=1,2$ there is a path
$\gamma_l$ in $\Gamma$ from $h_l$ to $g_l$ such that the label $U_l$
of $\gamma_l$ is a strongly Artin-reduced word and such that the paths
$\gamma_l$ and $\alpha_l$ are 4-close in $X$, as shown in
Figure~\ref{Fi:main}.

If $d_\Gamma(h_1,h_2)\le 1$ then $a$ is either trivial or a single
letter of $A\cup A^{-1}$. Let $W_1$ be the freely reduced form of the
word $aU_1$. Since $U_1$ is strongly Artin-reduced, the word $W_1$ is
Artin-reduced.

Assume now that $d_\Gamma(h_1,h_2)\ge 2$, so that $a\in G_{ij}$ for
some $m_{ij}<\infty$. If $a$ belongs to the cyclic subgroup generated
by one of the letters of $A$ (in which case $G_{ij}$ may not be
uniquely determined by $a$), say $a=a_i^q$, we let $W_1$ be the freely
reduced from of the word $a_i^q U_1$. Again, the word $W_1$ is
Artin-reduced because $U_1$ is strongly Artin-reduced.  Assume now
that $a$ is not a power of the generator, so that $a\in G_{ij}$ for a
unique pair $\{i,j\}$.  Let $U_1=Y_1 Z_1$ where $Y_1$ is the maximal
initial segment of $U_1$ which is a word in $G_{ij}$. Let $Y_1'$ be
the reduced word of smallest possible syllable length in $G_{ij}$
representing the element $aY_1$. Then $Y_1'$ and $Z_1$ are strongly
Artin-reduced and the word $W_1=Y_1'Z_1$ is Artin-reduced.

In either of the above cases, using the cone-vertex $v(h_1G_{ij})$ it
is easy to see that the path $\gamma_1$ is 2-close in $X$ to the path
$\lambda_1$ in $\Gamma$ from $h_2$ to $g_1$ labeled by $W_1$.

Similarly, by considering the product $U_2b$ we can find an
Artin-reduced word $W_2$ and a path $\lambda_2$ from $h_2$ to $g_1$ in
$\Gamma$ with label $W_2$ such that $\gamma_2$ and $\lambda_2$ are
2-close in $X$.

Now consider an equality diagram for the equality $W_1=_G W_2$ of
smallest area over the infinite presentation \eqref{art1} of $G$.
Applying Lemma~\ref{shape} to this diagram and using cone vertices we
see that $\lambda_1$ and $\lambda_2$ are 2-close in $X$. Adding the
distances between the paths considered in Figure~\ref{Fi:main} we see
that $\alpha_1$ and $\alpha_2$ are 14-close in $X$ and the Claim is
verified.
\end{proof}

\begin{rem}
  The value $K = 14$ obtained in the above proof of Theorem~\ref{A}
  using Papasolglu's criterion does not depend on the particular group
  $G$.  Thus all the coned-off Cayley graphs $X$ of the groups $G$
  covered by Theorem A are all $\delta$-hyperbolic for some $\delta>0$
  which is independent of $X$ and hence of the choice of $G$.  So the
  relative hyperbolicity of $G$ with respect to $G_{ij}$'s is in a
  sense uniform over the class of Artin groups covered by the Theorem.
\end{rem}

 \end{document}